\documentclass[reqno,11pt]{amsart}

\usepackage{epsf}
\usepackage{graphics}
\usepackage{graphicx}
\usepackage{amssymb}
\usepackage{amsmath}
\usepackage{mathcomp}
\usepackage{wasysym}
\date{}

\theoremstyle{plain}
\newtheorem{theorem}{Theorem}

\theoremstyle{definition}

\theoremstyle{remark}

\newtheorem*{remark}{Remark}

\def\F{{\mathbb F}}
\def\N{{\mathbb N}}

\def\Q{{\mathbb Q}}
\def\R{{\mathbb R}}

\title[A filtration of the Gordian graph] {Symmetric quotients of knot groups and a filtration of the Gordian graph}

\author{Sebastian Baader \and Alexandra Kjuchukova}

\begin{document}

\begin{abstract} We define a metric filtration of the Gordian graph by an infinite family of 1-dense subgraphs. The n-th subgraph of this family is generated by all knots whose fundamental groups surject to a symmetric group with parameter at least n, where all meridians are mapped to transpositions. Incidentally, we verify the Meridional Rank Conjecture for a family of knots with unknotting number one yet arbitrarily high bridge number.
\end{abstract}
\maketitle
\section{Introduction}

The Gordian graph $G$ is a countable graph whose vertices correspond to smooth knot types and whose edges correspond to pairs of knots related by a crossing change in a suitable diagram. The set of vertices of the Gordian graph carries a natural metric induced by the minimal path length between vertices. This metric, called the Gordian metric, generalizes the classical unknotting number, which is defined as the minimal number of crossing changes needed to transform a knot $K$ into the trivial knot $\Circle$, that is, $u(K)=d_G(K, \Circle)$ \cite{Mu}. The Gordian graph $G$ is locally infinite: for any given vertex of $G$, i.e. for any knot type, we obtain an infinite number of neighboring vertices by taking the connected sum of that knot with all knots whose unknotting number is one. Furthermore, every vertex of $G$ is contained in an arbitrarily large complete subgraph of $G$~\cite{HU}. In contrast, little is known about the global structure of the Gordian graph, except for the fact that it contains lattices of arbitrarily high rank~\cite{GG}. The purpose of this note is to construct an infinite descending sequence of subgraphs $G_n \subset G$ that captures the global geometry of the Gordian graph, by looking at finite symmetric quotients of the fundamental group of knot exteriors.

Fix a natural number $n \geq 2$. We say that a knot $K$ is ${n} \choose {2}$-colorable, if there exists a surjective homomorphism from the knot group $\pi_1(\R^3 \setminus K)$ onto the symmetric group $S_n$, mapping meridians to transpositions. Let $G_n \subset G$ to be the induced subgraph whose vertices correspond to knot types $K$ which are ${m} \choose {2}$-colorable, for some $m \geq n$. By definition, the sequence of subgraphs $G_n$ forms a descending chain, in fact strictly descending for $n\geq 2$, as we will see. A chain of graphs $\Gamma=\Gamma_1 \supset \Gamma_2 \supset \Gamma_3 \cdots$ forms a 1-dense metric filtration of the graph $\Gamma$, if for all $n\in \N$ the following three statements hold: 
\\
\begin{enumerate}
\item every vertex of $\Gamma$ is connected to a vertex of $\Gamma_n$ by an edge,
\item the maps $\Gamma_n \hookrightarrow \Gamma$ are isometric inclusions,
\item $\underset{m \in \N} \cap \Gamma_m=\emptyset$.
\end{enumerate}

The existence of a 1-dense metric filtration for a graph $\Gamma$ implies that every pair of non-neighboring vertices is connected by infinitely many different shortest paths. In particular, $\Gamma$ cannot be disconnected by removing finitely many vertices. The Gordian graph $G$ was already known to have these features~\cite{Ba}.

\begin{theorem} The chain of subgraphs $G=G_1 \supset G_2 \supset G_3 \cdots$ defined above forms a 1-dense metric filtration of the Gordian graph $G$.
\end{theorem} 

Thanks to properties (1) and (2) of a 1-dense metric filtration, determining the Gordian distance on $G$ is equivalent to determining its restriction to any of the subgraphs $G_n$, up to an error of two. Distance in $G_n$ can be studied via irregular branched covers of the 3-sphere along knots. In the first interesting case, $n=3$, these covers have been used to define a very effective knot invariant, as follows. To a surjective homomorphism  $\varphi: \pi_1(\R^3 \setminus K)\twoheadrightarrow S_3$ corresponds a three-fold irregular covering space $M$ of $S^3$ branched along $K$. When the two connected lifts of $K$ to $M$ represent torsion elements in homology, their linking number is a well-defined rational number, determined by $\varphi$. The set $lk(K)\subset\mathbb{Q}$ of linking numbers associated to all surjective homomorphisms of $\pi_1(\R^3 \setminus K)$ onto the symmetric group $S_3$ is called the {linking number invariant} of $K$~\cite{Re, Pe2}. We hope to obtain lower bounds on the Gordian distance of knots by estimating the effect of crossing changes on $lk(K)$. An outline of this new method is contained in the final section of this note.

\section{Symmetric quotients of knot groups}

The fundamental group of a knot complement admits a finite presentation via generators and relations, namely the Wirtinger presentation obtained from a knot diagram~\cite{CF}. This allows for a simple algorithmic decision of whether a given finite group $H$ is a quotient of $\pi_1(\R^3 \setminus K)$ or not. Consider the case $H=S_3$. In order to define a surjection onto $S_3$, the meridians of a knot need to be mapped to the three transpositions $(12),(13),(23)$, commonly referred to as colors. The group relations thus translate into the famous Fox 3-coloring conditions at the crossings of a diagram~\cite{Fo}. More generally, surjections onto dihedral groups $D_p$ correspond to Fox p-colorings of diagrams respecting similar rules.

The natural isomorphism $S_3 \simeq D_3$ gives rise to an alternative generalization of Fox 3-colorings, namely ${n} \choose {2}$-colorings, introduced earlier. 
Such a coloring encodes an ${n} \choose {2}$-representation of the knot group, that is, a surjective homomorphism
$$\pi_1(\R^3 \setminus K) \twoheadrightarrow S_n,$$
with the additional assumption that meridians be mapped to transpositions. We define the permutation number of a knot $K$ as follows:
$$p(K)=\max \{n \in \N \, | \, \pi_1(\R^3 \setminus K) \text{ admits an $\displaystyle{{n}\choose{2}}$-representation} \}.$$
We observe that the trefoil knot $3_1$ admits a Fox 3-coloring, in other words a ${3}\choose{2}$-representation, but no higher order ${n}\choose{2}$-represen\-tation, since the group $\pi_1(\R^3 \setminus 3_1) $ is generated by two meridians, whereas $S_n$ is not generated by two transpositions for $n \geq 4$. Therefore, $p(3_1)=3$. Similarly, $p(4_1)=2$ for the figure-eight knot $4_1$, since it does not admit a non-trivial 3-coloring, and its fundamental group is generated by two meridians. These simple observations extend to obtain a general upper bound on the permutation number of knots, from the minimal number of meridians needed to generate the knot group. Let $b(K)$ be the minimal bridge number of a knot $K$, defined as the minimal number of local maxima of the height function  among all representatives of $K \subset \R^3$. The fundamental group of a knot $K$ is generated by $b(K)$ meridians. Now the symmetric group $S_n$ cannot be generated by  fewer than $n-1$ transpositions. Indeed, in order to act transitively on $n$ numbers, one needs at least $n-1$ transpositions. Therefore, we obtain the following upper bound on the permutation number $p(K)$:
$$p(K) \leq b(K)+1.$$

Going back to the definition of the chain of subgraphs $G_n$ of the Gordian graph, we conclude that $K \notin G_n$ for $n > b(K)+1$, which implies property (3) of a metric filtration:
$$\underset{n \in \N} \cap G_n=\emptyset.$$

Our next goal is to construct knots with arbitrarily high permutation number yet unknotting number one. As a consequence, we will prove property (1) of a metric filtration, i.e. the 1-density of all subgraphs $G_n \subset G$. Fix a natural number $n$ and let $K_n$ be a knot with $p(K_n) \geq n$ and $u(K_n)=1$. Then the connected sum of knots $K \# K_n$ of any knot $K$ with $K_n$ is contained in $G_n$, since we can extend the existing homomorphism $\pi_1(\R^3 \setminus K_n) \twoheadrightarrow S_n$ to $\pi_1(\R^3 \setminus K \# K_n)$ by mapping all the meridians of $K$ to the transposition associated with the meridian of $K_n$ to which $K$ is attached. Moreover, the knot $K \# K_n$ is related to $K$ by a single crossing change, since the knot $K_n$ has unknotting number one. In other words, every knot (vertex) in $G$ is connected by an edge to a knot (vertex) in $G_n$. We are left to construct, for each $n\geq 3$, a knot $K_n$ with $p(K_n) \geq n$ and $u(K_n)=1$.

We first construct a family of knots with increasing permutation number. In fact, the $n$-times iterated connected sum of trefoil knots $3_1^n=3_1 \# \cdots \# 3_1$ will do, as $p(3_1^n) = n+2$. This can be seen by representing the knot $3_1^n$ as the closure of the braid $\sigma_1^3 \sigma_2^3 \cdots \sigma_{n-1}^3$ in the braid group $B_n$. Mapping the meridians around the $n$ bottom strands of that braid to the transpositions $(12),(13),\ldots,(1n)$, in this order, extends to a surjective homomorphism $\pi_1(\R^3 \setminus 3_1^n) \twoheadrightarrow S_{n+2}.$ The case $n=2$ is shown in Figure~1. This shows that $p(3_1^n) \geq n+2$. On the other hand, since $b(3_1^n)=n+1,$ by the previous discussion we  have  $p(3_1^n) \leq n+2$. In particular, $3_1^n\in G_{n+2}\setminus G_{n+3}$, establishing $G_{n}\supsetneq G_{n+1}$ for $n\geq 3$.

\begin{figure}[ht]
\qquad \qquad \scalebox{1.2}{\raisebox{-0pt}{$\vcenter{\hbox{\includegraphics{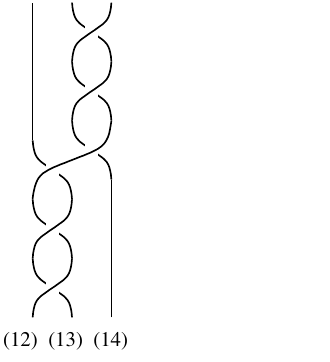}}}$}}
\caption{A ${{4}\choose{2}}$-coloring of the knot $3_1 \# 3_1$}
\end{figure}

In order to construct knots with unknotting number one and ${n} \choose {2}$-represen\-tations for $n$ large, we take suitable Whitehead doubles of iterated connected sums of the torus knot with 5 crossings. Let $T(2,5)$ be the torus knot represented as the closure of the braid $\sigma_1^5\in B_2$. Mapping the meridians of the two bottom strands of this braid to the two 3-cycles $(123),(345) \in A_5$ extends to a unique homomorphism $\pi_1(\R^3 \setminus T(2,5)) \rightarrow A_5$. Here the orientation of the meridians matters; we choose the convention that meridians cross under braid strands from the right to the left, as indicated by the arrow in Figure~2. Now let $K_1$ be the twisted Whitehead double of $T(2,5)$, defined as the closure of the braid
$$\beta_1=(\sigma_2 \sigma_1 \sigma_3 \sigma_2)^5 \sigma_3 \sigma_1 \in B_4,$$
with an additional clasp, as depicted on the right of Figure~2. The same figure exhibits a homomorphism $\pi_1(\R^3 \setminus K_1) \twoheadrightarrow S_5$, which maps the meridians of the bottom 4 strands of $\beta_1$ to the transpositions (12),(23),(34),(45). This doubling construction works since each 3-cycle can be written as a product of two transpositions.
\begin{figure}[ht]
\qquad \qquad \qquad \scalebox{1.0}{\raisebox{-0pt}{$\vcenter{\hbox{\includegraphics{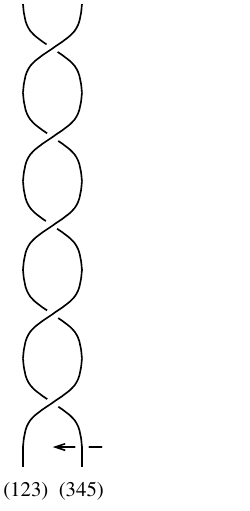}}}$}}
\quad
\scalebox{1.2}{\raisebox{-0pt}{$\vcenter{\hbox{\includegraphics{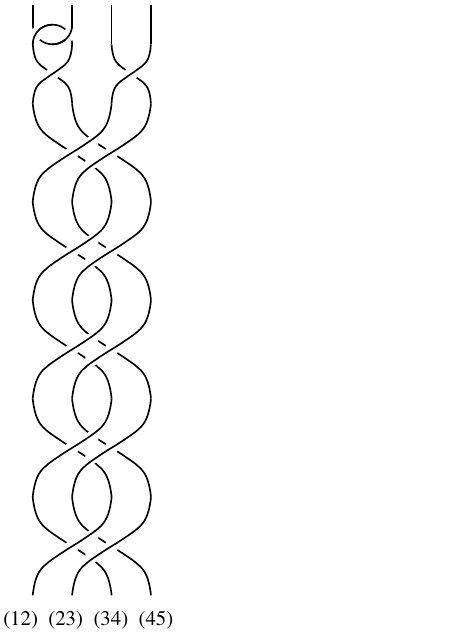}}}$}}
\caption{Colorings of $T(2,5)$ and its Whitehead double}
\end{figure}

The above argument generalizes to provide a homomorphism $$\pi_1(\R^3 \setminus T(2,5)^m) \rightarrow A_{2m+3}$$ mapping the meridian about the $k$-th bottom strand of the braid $\sigma_1^5\ldots \sigma_m^5\in B_{m+1}$ to the cycle $(2k-1~2k~2k+1)$. We define $K_m$ to be the Whitehead double of the knot $T(2,5)^m$, represented as the closure of the following braid in $B_{2m+2}$:
$$\beta_m=(\sigma_2 \sigma_1 \sigma_3 \sigma_2)^5 \sigma_3  (\sigma_4 \sigma_3 \sigma_5 \sigma_4)^5 \sigma_5 \ldots (\sigma_{2m} \sigma_{2m-1} \sigma_{2m+1} \sigma_{2m})^5 \sigma_{2m+1} \sigma_1,$$
again with an additional clasp on top of the first two strands.
By construction, the sequence of knots  $K_m$ have unknotting number one and are contained in $G_n \subset G$, for all $n \leq 2m+3$, since they admit an ${2m+3} \choose {2}$-coloring. This concludes the proof of the 1-density for all subgraphs $G_n \subset G$.

\begin{remark} The inequality $p(K) \leq b(K)+1$ between the permutation number and the bridge number of a knot implies that the bridge number of the knot $K_m$ defined above is at least $2m+2$. This bound is in fact sharp, since these knots are represented as closures of braids with $2m+2$ strands, i.e. by diagrams with precisely $2m+2$ local maxima. Moreover, we note that the permutation number provides the same lower bound for the meridional rank $\mu(K)$ of a knot, defined as the minimal number of meridians needed to generate the knot group: $$p(K)\leq\mu(K)+1.$$ Since $\mu(K)\leq b(K)$, this implies equality between the bridge number and the meridional rank of the knots $K_m$, and settles the Meridional Rank Conjecture for this family of knots (see Problem~1.11 in Kirby's list~\cite{Ki}). The latter conclusion could not have been drawn by using an analogous argument with Fox $p$-colorings. The mere existence of a non-trivial $p$-coloring does not provide an effective upper bound for the bridge number, due to the fact that every dihedral group is generated by two reflections. 
Moreover, counting $p$-colorings does not help either, since Whitehead doubles do not admit multiple independent non-trivial $p$-colorings. 
\end{remark}

\section{Constructing colored shortest paths}

In this section we prove that the inclusion maps $G_n \hookrightarrow G$ are isometric inclusions, for all $n \in \N$. In other words, the intrinsic path metric on $G_n$ coincides with the path metric on the ambient space $G$. The strategy of proof is as follows: let $A,B \in G_n$ be two knots which admit ${a} \choose {2}$- and ${b} \choose {2}$-colorings ($a,b \geq n$), respectively, and let $d_G(A,B)=m \in \N$ be the Gordian distance between them. Then there exists a sequence of knots $K_0,K_1,\ldots K_m$ successively related by crossing changes, i.e. $d_G(K_i,K_{i+1})=1$ for all $i \leq m-1$, with $K_0=A$, $K_m=B$. We will construct a knot $\widetilde{K}_1 \in G_n$ with $d_G(K_0,\widetilde{K}_1)=d_G(\widetilde{K}_1,K_2)=1$. By repeating this argument inductively, we will end up with a path from $A$ to $B$ in $G_n$ of the same length as the original one.

In order to construct the knot $\widetilde{K}_1$, we will make use of the fact that the two crossing changes relating $K_1$ to each of $K_0$ and $K_2$ can be realized in the same diagram of $K_1$, as explained in~\cite{Ba}. Here is an outline of the argument: a crossing change between two strands can be represented by a framed chord with endpoints on the knot. The endpoints of these chords can be moved along the knot, and moreover they can be contracted to small segments. We may therefore assume that $K_1$ has a diagram with a section as depicted in Figure~3, containing two neighboring clasps, a crossing change at which transforms $K_1$ into $K_0$ and $K_2$, respectively (see~\cite{Ba} for details).
\begin{figure}[ht]
$K_1$ \, 
\scalebox{1.2}{\raisebox{-0pt}{$\vcenter{\hbox{\includegraphics{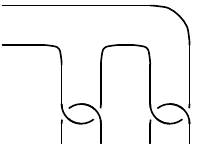}}}$}}

\bigskip
$K_0$ \, 
\scalebox{1.2}{\raisebox{-0pt}{$\vcenter{\hbox{\includegraphics{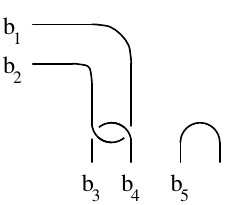}}}$}}
\quad $K_2$ \, 
\scalebox{1.2}{\raisebox{-0pt}{$\vcenter{\hbox{\includegraphics{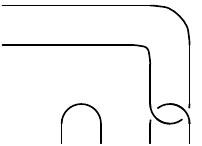}}}$}}
\caption{Diagrams of the knots $K_1,K_0,K_2$}
\end{figure}

Our assumption on $A=K_0$ tells us that the diagram of $K_0$ admits a coloring by transpositions that generate the symmetric group $S_a$. The section of the diagram on the bottom left of Figure~3 has five connected arcs $b_1,b_2,b_3,b_4,b_5$, whose meridians are sent to five transpositions, some of which may coincide. Depending on these transpositions, we will define a new knot $\widetilde{K}_1$, which admits an ${a} \choose {2}$-representation, as well. We distinguish two cases:

\smallskip
\noindent
1) The bridges $b_1$, $b_2$ are sent to the same transposition $(ij)$. Then the Wirtinger relations at the clasp imply that $b_3$, $b_4$ are also sent to one transposition $(kl)$, which commutes with $(ij)$. In this case, we define the knot $\widetilde{K}_1$ as in Figure~4, where it is assumed that $\widetilde{K}_1$ coincides with $K_1$ outside the depicted region. Observe that the ${a} \choose {2}$-coloring of $K_0$ carries over to an ${a} \choose {2}$-coloring of $\widetilde{K}_1$. Therefore, $\widetilde{K}_1$ is in $G_a$, in turn in $G_n$, since $a \geq n$. Moreover, this new knot is still related to each of $K_0$ and $K_2$ by a single crossing change.
\begin{figure}[ht]
\scalebox{1.4}{\raisebox{-0pt}{$\vcenter{\hbox{\includegraphics{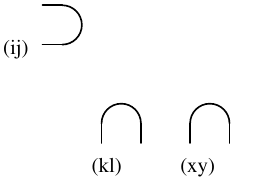}}}$}}
\caption{$\widetilde{K}_1$ in case 1}
\end{figure}

\smallskip
\noindent
2) The arcs $b_1$, $b_2$ are sent to different transpositions, whose support overlaps in one number, $(jk)$ and $(ij)$, since $b_1$ and $b_2$ meet at a crossing. Then the Wirtinger relations at the clasp imply that $b_3$, $b_4$ are sent to $(jk)$ and $(ik)$. We distinguish three subcases, depending on the transposition $(xy)$ associated with the arc $b_5$.

\smallskip
\noindent
a) The transposition $(xy)$ commutes with $(jk)$, that is, $(xy)=(jk)$ or $\{x,y\}\cap\{j,k\}=\emptyset$.  
In this case, the ${a} \choose {2}$-coloring of $K_0$ carries over to an ${a} \choose {2}$-coloring of $K_1$, so we keep $\widetilde{K}_1=K_1$.

\smallskip
\noindent
b) The transposition $(xy)$ commutes with $(ij)$, that is, $(xy)=(ij)$ or $\{x,y\}\cap\{i,j\}=\emptyset$. We define the knot $\widetilde{K}_1$ as in Figure~5. By construction, the ${a} \choose {2}$-coloring of $K_0$ carries over to an ${a} \choose {2}$-coloring of $\widetilde{K}_1$. Moreover, $\widetilde{K}_1$ is related to each of $K_0$ and $K_2$ by a single crossing change.
\begin{figure}[ht]
\scalebox{1.4}{\raisebox{-0pt}{$\vcenter{\hbox{\includegraphics{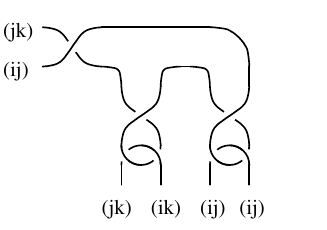}}}$}}
\caption{$\widetilde{K}_1$ in case 2b, with $(xy)=(ij)$}
\end{figure}

\smallskip
\noindent
c) The transposition $(xy)$ commutes with $(ik)$, that is, $(xy)=(ik)$ or $\{x,y\}\cap\{i,k\}=\emptyset$. We define the knot $\widetilde{K}_1$ as in Figure~6. As before, the ${a} \choose {2}$-coloring of $K_0$ carries over to an ${a} \choose {2}$-coloring of $\widetilde{K}_1$ and $\widetilde{K}_1$ is related to each of $K_0$ and $K_2$ by a single crossing change.
\begin{figure}[ht]
\scalebox{1.4}{\raisebox{-0pt}{$\vcenter{\hbox{\includegraphics{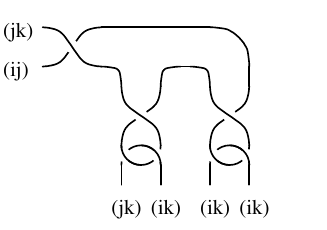}}}$}}
\caption{$\widetilde{K}_1$ in case 2c, with $(xy)=(ik)$}
\end{figure}

\smallskip
\noindent
Note that the above cases exhaust all possibilities, with some redundancy. For example, $(xy)$ with $\{x,y\}\cap\{i, j, k\}=\emptyset$ is covered by 2 a), b), c). Where there is overlap, one can choose which construction to use. This concludes the proof of Theorem~1.

\section{Speculation on the Gordian distance}

In this final section we restrict our attention to 3-colorable knots. We say a knot is twice three-colorable if it admits two independent non-trivial 3-colorings, i.e. 3-colorings that are not related by a permutation of the three transpositions $(12),(13),(23)$. Examples of twice 3-colorable knots are the square knot $3_1 \# \overline{3}_1$, the granny knots $3_1 \# 3_1$ and $\overline{3}_1 \# \overline{3}_1$, and the pretzel knot $P(3,3,3)$. We will outline an obstruction for a knot to be related by a crossing change to a twice 3-colorable knot. As an application, we will show that the square knot and the knot $8_{21}$ are not related by a crossing change. This is a hard case, since the knot $8_{21}$ is related to the trivial knot by a single crossing change, and the square knot is indistinguishable from the trivial knot by the well-known concordance knot invariants $\sigma,s,\tau$, because it is a ribbon knot.

A well-known features of the 3-colorings of a knot is that they from an $\F_3$-vector space, after any identification of the three transpositions $(12)$, $(13)$, $(23)$ with the numbers $0,1,2 \in \F_3$. This means that the coloring conditions at crossings are preserved under addition, a fact that is easily checked. Now suppose that a knot $K$ has two independent non-trivial 3-colorings $C_1,C_2$. Then the sum $C_1+C_2$ and the difference $C_1-C_2$ are also non-trivial 3-colorings of $K$. Moreover, given any crossing $X$ in a diagram of $K$, one of the four 3-colorings $C_1,C_2,C_1+C_2,C_1-C_2$ will have two coinciding colors at $X$, hence it will be monochromatic at $X$. As a consequence, every knot $K_2$ that is related to a twice 3-colorable knot $K_1$ by a single crossing change will be related to $K_1$ by a monochromatic crossing change with respect to suitable non-trivial 3-colorings of $K_1,K_2$. In this case, the 3-fold irregular branched coverings associated with the representations $\pi_1(\R^3 \setminus K_{1,2}) \twoheadrightarrow S^3$ will be related by controlled surgeries. We refer the reader to~\cite{Pe1} for a precise definition of these branched covering spaces. Whenever the branch curve lifts to a curve of finite order in homology, the linking number of the pair of lifts of the branch curve is a well-defined rational number. The set of linking numbers obtained by considering all non-trivial 3-colorings of a knot $K$ forms the knot invariant $\text{lk}(K) \subset \Q$ mentioned in the introduction. 

We conjecture that monochromatic crossing changes do not affect linking numbers in branched coverings. This was recently confirmed by Perko~\cite{Pe3} for a special case of 3-colored knot diagrams. If true in general, this conjecture would imply that whenever a knot $K_1$ is related to a twice 3-colorable knot $K_2$ by a single crossing change, then the linking number invariants $\text{lk}(K_1)$ and $\text{lk}(K_2)$ must have at least one number in common. Using the algorithm implemented in~\cite{CK}, we computed a list of linking numbers for the knots $3_1,8_{20},8_{21}$ and their mirror images, as well as for the twice 3-colorable granny and square knots:
$$\text{lk}(3_1)=\{2\}, \, \text{lk}(\overline{3}_1)=\{-2\}, \, \text{lk}(8_{20})=\text{lk}(\overline{8}_{20})=\{0\},$$
$$\text{lk}(8_{21})=\{4\}, \, \text{lk}(\overline{8}_{21})=\{-4\},$$
$$\text{lk}(3_1 \# 3_1)=\{2,4\}, \, \text{lk}(\overline{3}_1 \# \overline{3}_1)=\{-4,-2\}, \, \text{lk}(3_1 \# \overline{3}_1)=\{-2,0,2\}.$$
We found crossing changes relating knots from the first two lines to knots from the third line in all the cases that have a linking number in common. In all other cases, the obstruction described above would rule out a single crossing change relating a knot from the first two lines to a knot from the third line.

\bigskip
\noindent
Mathematisches Institut, Universit\"at Bern, Sidlerstrasse 5, CH-3012 Bern, Switzerland

\smallskip
\noindent
\texttt{sebastian.baader@math.unibe.ch}

\bigskip
\noindent
Mathematics Department, University of Wisconsin--Madison, 480 Lincoln Dr, Madison, WI 53703

\smallskip
\noindent
\texttt{kjuchukova@wisc.edu}

\end{document}